\documentclass{amsart}
\usepackage{amsmath,amscd,amssymb,amsthm,epsfig,psfrag,latexsym,subfigure}

\newtheorem{theorem}{Theorem}[section]

\newcommand{\ZZ}{\mathbb{Z}}
\newcommand{\DD}{\mathbb{D}}
\newcommand{\NN}{\mathbb{N}}
\newcommand{\QQ}{\mathbb{Q}}
\newcommand{\CC}{\mathbb{C}}
\newcommand{\RR}{\mathbb{R}}

\def\vol{\mbox{\rm{Vol}}}

\def\a{\alpha}
\def\b{\beta}
\def\d{\partial}

\def\g{\gamma}
\def\i{\iota}

\def\G{\Gamma}
\def\th{\theta}
\def\l{\lambda}

\def\i{\iota}

\def\S{\Sigma}

\def\tr{\mbox{\rm{tr}}}

\def\PSL{\mbox{\rm{PSL}}}

\def\SL{\mbox{\rm{SL}}}

\def\min{\mbox{\rm{min}}}

\def\Int{\mbox{int}\ }

\def\HH{\mathbb{H}}
\def\AA{\mathbb{A}}

\def\MC{\mathcal{C}}

\def\CC{\mathbb{C}}

\edef\t@mp{\catcode`\noexpand\#=\the\catcode`\#}%
    \catcode`\#=12
    \def\h@sh{#}%
    \t@mp

\edef\t@mp{\catcode`\noexpand\~=\the\catcode`\~}%
    \catcode`\~=12
    \def\tild@{~}%
    \t@mp

\begin{document}
\Large
\title{Transcendental ending laminations}
\author{Ian Agol}
\address{Department of Mathematics, University of Illinois at
Chicago, 322 SEO m/c 249, 851 S. Morgan Street, Chicago, IL
60607-7045}



\email{agol@math.uic.edu}
\thanks{Partially supported by NSF grant DMS 0204142 and the Sloan
Foundation}

\begin{abstract}
Yair Minsky showed that punctured torus groups are classified by a
pair of ending laminations $(\nu_-,\nu_+)$. In this note, we show
that there are ending laminations $\nu_+$ such that for any choice
of $\nu_-$, the punctured torus group is transcendental as a
subgroup of $\PSL_2\CC$.
\end{abstract}

\maketitle

\section{Introduction}
We investigate in this note the algebraic properties of Kleinian
groups, that is finitely generated, discrete subgroups of
$\PSL_2\CC$. Let $\AA\subset \CC$ be the algebraic numbers. A
subgroup $G<\PSL_2\CC$ is {\it algebraic} if it can be conjugated
into $\PSL_2\AA$. Otherwise, it is {\it transcendental}. By Mostow
rigidity, all finite covolume Kleinian groups are algebraic. In
fact, it follows from the proof of Thurston's hyperbolization
theorem for Haken manifolds that any Kleinian group has {\it some}
faithful representation into $\PSL_2\CC$ whose image is algebraic.
This is proven by embedding the abstract Kleinian group into the
fundamental group of a finite volume hyperbolic 3-orbifold. The
induced hyperbolic structure on the subgroup from Thurston's
construction is geometrically finite \cite{Mo}. It seems possible
that there are algebraic Kleinian groups which are geometrically
finite, but which are not the subgroup of any finite volume
Kleinian group. However, it follows from the geometric tameness
theorem \cite{Ag04, CG04} and the covering theorem \cite{Ca} that
if a non-geometrically finite Kleinian group $\G$ is a subgroup of
a finite volume Kleinian group $\G'$, then $\G$ is a surface group
containing a  fiber subgroup $\G_0$ of a virtual fibration of
$\HH^3/\G'$ of index at most 2 in $\G$. Each end of a Kleinian
group is either geometrically finite, or simply degenerate
\cite{Ca93}. Thus, in this case $\G_0$ has two simply degenerate
ends, and there is an infinite order isometry $\psi\in \G'-\G_0$
which normalizes $\G_0$, and acts by a translation on
$\HH^3/\G_0$. It does not appear to be known whether there are
degenerate algebraic Kleinian groups which are not of this type. A
similar argument to the covering theorem shows that a singly
degenerate hyperbolic 3-manifold cannot immerse totally
geodesically in a higher dimensional geometrically finite
hyperbolic manifold. Any algebraic Kleinian group $\G$ embeds in a
natural way into a higher rank S-arithmetic lattice. But the
embedding does not induce a quasi-isometric embedding on the
convex core of $\HH^3/\G$ into the corresponding cover of the
symmetric space, otherwise the techniques of the covering theorem
would carry over. One may construct simplicial ruled surfaces in
higher rank symmetric spaces, but the bounded diameter lemma and
local finiteness fail.

It would be interesting to know if there can be singly degenerate
algebraic surface Kleinian groups. We demonstrate in this note
that there are families of degenerate punctured torus Kleinian
groups satisfying certain geometric restrictions which must be
transcendental. The families of groups we consider are
uncountable, and therefore all but countably many groups in a
family will be transcendental; the point of the result is that
{\it all} of the groups in the family are transcendental.

One motivation for studying this question is to understand whether
there is an algorithm to tell that a collection of matrices in
$\SL_2\AA$ generates a discrete subgroup of $\SL_2\CC$. If every
non-geometrically finite Kleinian group were the subgroup of a
finite volume Kleinian group, then an algorithm would exist. So
the discovery of a singly degenerate algebraic Kleinian group
would demonstrate the difficulty of this question.

{\bf Acknowledgment} We thank Caltech and Nathan Dunfield for
their hospitality during a visit in which this work was carried
out.
\section{Ending lamination theorem}
In this section, we review the classification of punctured torus
Kleinian groups from \cite{Minsky}.
 Let $\S$ be a punctured torus.  A punctured torus group
$G$ is a 2 generator free discrete subgroup of $\PSL_2\CC$, such
that there exist two generators of $G$ whose commutator is
parabolic. $G$ is naturally identified as the image of a discrete
faithful representation $\rho:\pi_1(\S)\to \PSL_2\CC$ such that
the monodromy of the puncture is parabolic. Then $N=\HH^3/G$ is a
punctured torus manifold. Bonahon showed that $N\cong \S\times
\RR$ \cite{Bon}.

\begin{figure}[htbp]
    \begin{center}
    \epsfbox{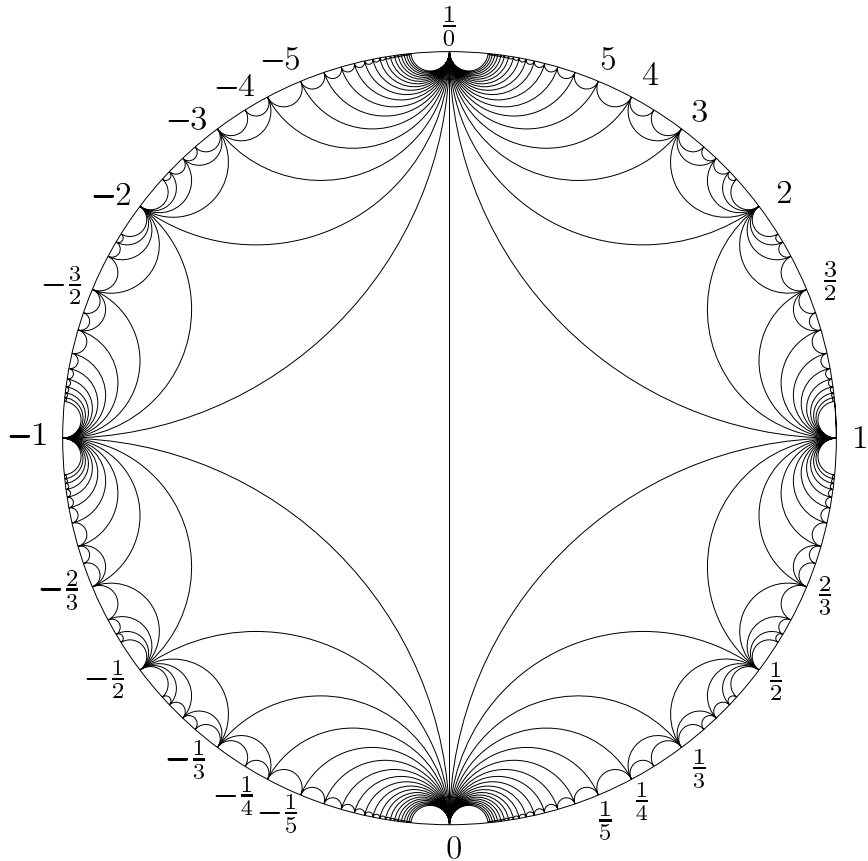}
    \caption{\label{Farey} The Farey graph in the unit disk,
    with some of the vertices labelled}
    \end{center}
\end{figure}

Let $\HH^2$ denote the upper half plane model for hyperbolic
space, with boundary $\RR$. The Farey triangulation is a graph in
$\overline{\HH^2}=\HH^2\cup \RR \cup \infty$ with vertices at
$\hat{\QQ}=\QQ\cup \infty$, and two rational numbers $\frac{p}{q},
\frac{r}{s}$ are connected by a hyperbolic geodesic if $|ps-qr|=1$
(where we let $\infty =\frac10$). The irrational numbers in $\RR$
represent the ending laminations. See figure \ref{Farey} for a
partial picture of the Farey triangulation in the disk model of
hyperbolic space.
The curve complex $\MC(\S)$ of the punctured torus has vertices
consisting of isotopy classes of simple closed curves. The edges
of $\MC(\S)$ join pairs of curves which intersect exactly once.
There is a natural bijection of $\MC(\S)^{(1)}$ with the Farey
triangulation, for example by considering projective homology
classes in $H_1(\S)\cong\ZZ^2$. Given a basis $\{\a,\b\}$ for
$H_1(\S)$, any oriented simple closed curve represents a homology
class $p\a+q\b$, with $\gcd(p,q)=1$. Then we get an element
$-\frac{p}{q}\in \hat{\QQ}$, which gives a bijection with
unoriented simple closed curves. The number $|pr-qs|$ represents
the geometric intersection number between two simple closed curves
$p/q, r/s$.

If we lift a punctured torus group representation $\rho:\pi_1\S\to
\PSL_2\CC$ to $\tilde{\rho}:\pi_1\S\to \SL_2\CC$, then the
function $\tr(\tilde{\rho}(\a))$ is well-defined on conjugacy
classes of elements in $\pi_1\S$. This function takes values in
$\AA$ if $\rho$ is algebraic. The lift $\tilde{\rho}$ is
well-defined up to multiplying generators by $-I$. Given three
adjacent curves $\a,\b,\g\in \MC(\S)$, with traces
$x=\tr(\tilde{\rho}(\a)), y=\tr(\tilde{\rho}(\b)),
z=\tr(\tilde{\rho}(\g))$, they satisfy the Markoff equation
\begin{equation}x^2+y^2+z^2=xyz.\end{equation} Moreover, if $\delta\in \MC(\S)$ is
adjacent to $\a,\b$ with trace $w=\tr(\tilde{\rho}(\delta))$, then
these satisfy the equation \begin{equation}w+z=xy.
\label{recursion}\end{equation} These formulae follow from trace
identities in $\SL_2\CC$, see \cite{Bowditch98} for details.

Let $\DD=\overline{\HH^2}$. Bonahon showed the existence of end
invariants $(\nu_-,\nu_+)\subset (\DD\times \DD)\backslash
\Delta$, where $\Delta$ is the diagonal of $\d\DD\times\d\DD$
\cite{Bon}. Teichm\"{u}ller space for the once punctured torus is
identified with $\Int \DD=\HH^2$. Minsky proved that the end
invariants of $G$ determine $G$ uniquely up to conjugacy in
$\PSL_2\CC$ (Theorem A, \cite{Minsky}). We may thus call $\nu_+$ a
{\it transcendental ending lamination} if the punctured torus
group associated to $(\nu_-,\nu_+)$ is transcendental for all
$\nu_-\in \DD-\{\nu_+\}$.

Let $N$ be a punctured torus manifold. If one removes the interior
of the Margulis tube associated to the parabolic commutator to get
$\hat{N}$, then $\hat{N}$ has two ends, $e_{\pm}$. The ends are
either geometrically finite, or simply degenerate. To each end,
associate an end invariant $\nu_{\pm}\in\DD$. If $e_s$ ($s\in
\{+,-\}$) is geometrically finite, then $\nu_s\in \Int \DD$ is the
point in Teichm\"{u}ller space of the conformal structure of the
domain of discontinuity corresponding to $e_s$, or if there is an
accidental parabolic, it is the curve $\nu_s\in \hat{\QQ}=\MC(\S)$
represented by the accidental parabolic in $e_s$. If $e_s$ is
simply degenerate, then there is a unique ending lamination
$\nu_s\in \RR-\QQ$ which is a limit point of simple curves in
$\MC(\S)$ whose geodesic representatives in $N$ have bounded
length and exit the end $e_s$.

\begin{figure}[htbp]
    \begin{center}
    \epsfbox{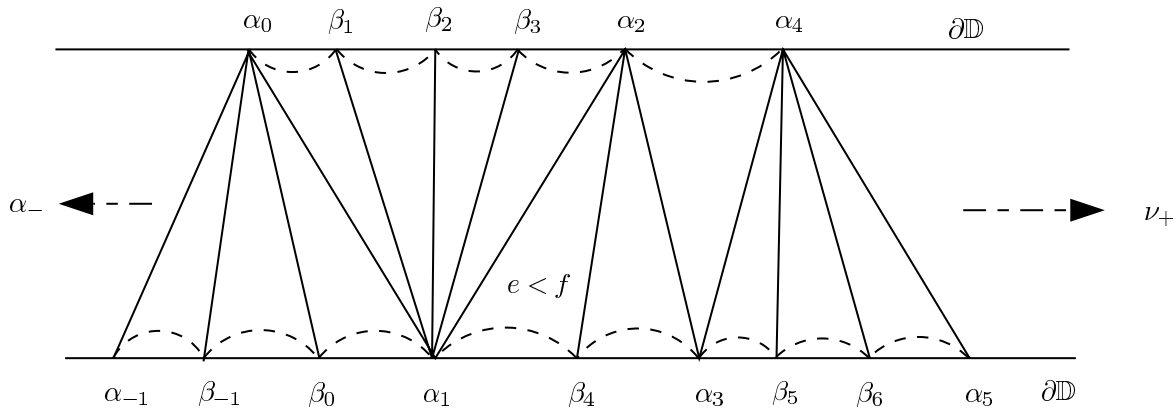}
    \caption{\label{pivot} The pivot
sequence }
    \end{center}
\end{figure}

To each end invariant $\nu_s$, we associate an element $\a_s$ of
$\hat{\RR}$, $s=\pm$. If the end $e_s$ is simply degenerate or has
an accidental parabolic, then $\a_s=\nu_s$. If $e_s$ is
geometrically finite, then $\a_s$ is a systole in the hyperbolic
structure on $\S$ corresponding to the Teichm\"uller point
$\nu_s$, which is defined up to finite ambiguity. We will assume
from now on that $e_+$ is simply degenerate, so that $\nu_+\in
\RR-\QQ$. Define $E=E(\a_-,\a_+)$ to be the set of edges of the
Farey graph which separate $\a_-$ from $\a_+$ in $\DD$. Let
$P_0\in \hat{\QQ}$ be the set of vertices of $\MC(\S)$ belonging
to at least two edges of $E$. The edges of $E$ admit a natural
order where $e<f$ if $e$ separates the interior of $f$ from
$\a_-$, and this induces an order on $P_0$ (see figure
\ref{pivot}). So we can arrange $P_0$ as a sequence
$\{\a_n\}_{n=\i}^{\infty}$, where $\i=-\infty$ if $\nu_-\in
\RR\backslash \QQ$, and $\i=0$ otherwise. We note that the
vertices of $E$ which are not elements of $P_0$ also are ordered,
and we will denote these by $\{\b_i\}_{i=\iota}^{\infty}$, with
the convention that $\b_1$ is adjacent to $\a_0$ in $\MC(\S)$ and
$\a_0<\b_1$ (meaning that any edge of $E$ incident with $\a_0$ is
$<$ the unique edge of $E$ incident with $\b_1$). The {\it width}
of $\a_i$ is $1+$ the number of vertices in $E\backslash P_0$
which are adjacent to $\a_i$ in $\MC(\S)$ and which are $>
\a_{i-1}$ and $< \a_{i+1}$ (in the example in figure \ref{pivot},
$w(0)=3, w(1)=4, w(2)=2, w(3)=1, w(4)=3$).

The pivot sequence $\{\a_n\}_{n=0}^{\infty}$ determines $\nu_+$
uniquely. In fact, given adjacent $\a_0,\a_1\in \MC(\S)$, and the
width sequence $\{w(n)\}_{n\in \NN}$, the pivot sequence
$\{\a_n\}$ and $\nu_+$ are uniquely determined. Assuming a fixed
lift of a punctured torus group $\tilde{\rho}:\pi_1\S\to
\SL_2\CC$, let us denote $a_i=\tr(\tilde{\rho}(\a_i)),
b_j=\tr(\tilde{\rho}(\b_j))$. From the trace relation equation
\ref{recursion}, we can compute $a_i$ or $b_j$ recursively as a
polynomial function of the traces of $\a_0,\a_1,\b_1$. In the
example in figure \ref{pivot}, $b_2=b_1a_1-a_0$,
$b_3=b_2a_1-b_1=(b_1a_1-a_0)a_1-b_1$.  We let $\l(\g)=l+i\th$ be
the complex translation length of an element (or conjugacy class)
$\g\in\PSL_2\CC$, normalized so that $l\geq 0, \th\in (-\pi,\pi]$.
It is determined by the identity $\tr^2\g=4\cosh^2(\l/2)$. Given
$\rho:\pi_1(\S)\to \PSL_2\CC$, we get a function on $\MC(\S)$
denoted $\l(\a)\equiv\l(\rho(\a))$.

\section{Construction of transcendental ending laminations}
 We will describe the construction of
width sequences $\{w(n)\}_{n=1}^{\infty}$ with associated ending
lamination $\nu_+$ such that for any choice of ending lamination
$\nu_-\in \DD-\{\nu_+\}$, the punctured torus group $G$ associated
to $(\nu_-,\nu_+)$ is transcendental. The construction is akin to
Cantor diagonalization and to Liouville's construction of
transcendental numbers. The non-trivial input which we take is a
result of Minsky.

\begin{theorem} (4.4,4.5 \cite{Minsky})
There exist constants $c_2, c_3, c_4$ such that
$$\frac{c_2}{w(n)^2}\leq l(\a_n)\leq \frac{c_3}{w(n)^2},$$
$$|w(n)-\frac{2\pi}{\th(\a_n)}|\leq c_4.$$
\end{theorem}

From this theorem, we clearly have $l(\a_n)=O(w(n)^{-2})$,
$\th(\a_n)=O(w(n)^{-1})$, so $\l(\a_n)=O(w(n)^{-1})$. It follows
that
\begin{equation}|a_n^2-4| =
|4\cosh^2(\l(\a_n)/2)-4|= O(w(n)^{-2}).\label{bound}\end{equation}

{\bf Construction:} Choose an exhaustion of $\AA$ by finite
subsets $\AA=\cup_{i=0}^{\infty} A_i$, $|A_i|<\infty$, which is
possible since $\AA$ is countable. Also, choose a sequence
$\{n_i\}\subset \NN$.  Assume by induction that
$\{w(1),...,w(n_i-1)\}$ have been defined. We may choose a triple
of generators $\a_0,\a_1,\b_1$ with traces $a_0,a_1,b_1$. Assume
that $a_0,a_1,b_1\subset A_i$. There is a polynomial
$P_i(x,y,z)\in \ZZ[x,y,z]$ such that $a_{n_i}=P_i(a_0,a_1,b_1)\neq
\pm 2$ (following from theorem 2.1, $\l(\a_n)\neq 0$ since any
parabolic element is conjugate to either the peripheral curve or
to an accidental parabolic corresponding to an end invariant by
theorem 4.1 \cite{Minsky}). Let
\begin{equation} m_i=\min\{ |P_i(x,y,z)^2-4|\ \ni\ x,y,z\in A_i, P_i(x,y,z)\neq
\pm 2\}.\end{equation}
 Choose $w(n_i)$ large enough that
$|a_{n_i}^2-4|< m_i,$ for any pivot  $\a_{n_i}$ with width
$w(n_i)$, which we may choose by equation \ref{bound}. We may
choose the widths $\{w(n_i+1),...,w(n_{i+1}-1)\}$ arbitrarily.
This describes the recursive construction of the sequence
$\{w(n)\}_{n=1}^{\infty}$.

\begin{theorem}
If $\nu_+$ is associated to the width sequence $\{w(1),w(2),...\}$
coming from the construction, then for any ending lamination
$\nu_-$, the punctured torus group $G$ associated to
$(\nu_-,\nu_+)$ is transcendental.
\end{theorem}

\begin{proof}
Suppose we had a punctured torus group $G < \PSL_2\AA$, with
associated ending pair $(\nu_-,\nu_+)$, where $\nu_+$ has the
width sequence $\{w(n)\}$ given in the construction above. Then
$\{a_0,a_1,b_1\}\subset A_i$, for some $i$. Thus, we have for the
pivot $\a_{n_i}$, $$0\neq |a_{n_i}^2-2| \geq m_i >
|a_{n_i}^2-4|,$$ a contradiction.
\end{proof}

\section{Conclusion}
Most likely, the numbers $\nu_+\in \RR-\QQ$ coming from the
construction are transcendental, since the continued fraction
coefficients of $\nu_+$ are given by $w(n)$ (up to ambiguity of a
finite initial sequence), and the numbers $w(n)$ should grow fast
enough to satisfy Liouville's criterion for being transcendental.
One could get explicit sequences $\{w(n)\}$ which are defined
recursively and satisfy the requirements of the construction by
using the generalized Liouville theorem, which bounds the length
of $P_i(a_0,a_1,b_1)$ in terms of the length of $P_i$ and the
lengths of $a_0,a_1,b_1$. The length of a polynomial is the sum of
the absolute values of its coefficients, and the length of an
algebraic number is the length of its minimal polynomial. Then the
set $A_i\subset \AA$ would be all algebraic numbers with length
and degree $\leq i$, which is a finite set, and one could estimate
the growth of the bound on $P_i(a_0,a_1,b_1)$ using equation
\ref{recursion} recursively. It might be interesting to compute
the asymptotics of $\{w(n)\}$ using this method.

The same technique produces transcendental ending laminations for
4 punctured sphere groups and 2 parabolic generator groups. It
would be interesting to try to use the techniques of the general
ending lamination theorem  to find transcendental ending
laminations for arbitrary Kleinian surface groups \cite{Minsky03,
BCM}.

\def\cprime{$'$} \def\cprime{$'$}

\end{document}